\documentclass[10pt]{amsart}
\usepackage{amsthm,amsmath,amssymb}

\global\overfullrule10pt

\theoremstyle{plain}
\newtheorem{theorem}{Theorem}
\newtheorem{corollary}[theorem]{Corollary}
\newtheorem{lemma}[theorem]{Lemma}
\newtheorem{proposition}[theorem]{Proposition}
\newtheorem*{theorem*}{Theorem}
\newtheorem*{corollary*}{Corollary}
\newtheorem*{proposition*}{Proposition}
\theoremstyle{definition}

\newtheorem{definition}[theorem]{Definition}
\theoremstyle{remark}
\newtheorem{remark}[theorem]{Remark}

\newcommand\CC{{\mathbf C}}
\newcommand\RR{{\mathbf R}}
\newcommand\ZZ{{\mathbf Z}}
\newcommand\NN{{\mathbf N}}

\newcommand\BB{{\mathbf B}}
\newcommand\spr[1]{\langle#1\rangle}
\newcommand\hol{{\text{hol}}}
\newcommand\harm{{\text{harm}}}

\renewcommand\Re{\operatorname{Re}}

\newcommand\oz{\overline z}
\newcommand\Bn{{\BB^n}}
\newcommand\pBn{{\partial\Bn}}
\newcommand\bn{{B^n}}
\newcommand\pbn{{\partial\bn}}
\newcommand\Hardy{{\text{\rm Hardy}}}

\newcommand\wt{\widetilde}
\newcommand\td{\wt\Delta}

\newcommand\FF[5]{{}_#1\!F_#2\Big(\begin{matrix}#3\\#4\end{matrix}\Big|#5\Big)}
\newcommand\Mh{$M$-harmonic }
\newcommand\Ph{pluriharmonic }
\newcommand\LL{\mathcal L}
\newcommand\chpq{\mathcal H^{pq}}
\newcommand\bhpq{\mathbf H^{pq}}
\newcommand\chp{\mathcal H^p}
\newcommand\bhp{\mathbf H^p}
\newcommand\Un{{U(n)}}
\newcommand\cdd[2]{#1\cdot\overline#2}
\newcommand\HH{{\mathcal H}}
\newcommand\dsph{\Delta_{\text{sph}}}
\newcommand\cR{\mathcal R}
\newcommand\cN{\mathcal N}

\makeatletter
\newcommand{\vast}{\bBigg@{4}}
\newcommand{\Vast}{\bBigg@{5}}
\makeatother

\renewcommand\[{\begin{equation}}
\renewcommand\]{\end{equation}}
\allowdisplaybreaks[4]

\newcommand\cP{\mathcal P}
\newcommand\cici{{\circ\circ}}
\renewcommand\AA{\mathcal A}
\newcommand\jedna{\mathbf 1}
\newcommand\itv[1]{[#1)}
\newcommand\autBn{\operatorname{Aut}(\Bn)}
\newcommand\MM{\mathcal M}

\newcommand\dbar{\overline\partial}
\newcommand\sq{{\square}}
\newcommand\XX{\mathcal X}
\newcommand\nct{\nabla_{\text{ct}}}

\begin{document}

\title[$M$-harmonic Dirichlet space]{The $M$-harmonic Dirichlet space on the ball}
\author[M.~Engli\v s]{Miroslav Engli\v s}
\address{Mathematics Institute, Silesian University in Opava,
 Na~Rybn\'\i\v cku~1, 74601~Opava, Czech Republic {\rm and }
 Mathematics Institute, \v Zitn\' a 25, 11567~Prague~1,
 Czech Republic}
\email{englis{@}math.cas.cz}
\author{El-Hassan Youssfi}
\address{Aix-Marseille Universit\'e, I2M UMR CNRS~7373,
 39 Rue F-Juliot-Curie, 13453~Marseille Cedex 13, France}
\email{el-hassan.youssfi{@}univ-amu.fr}
\thanks{Research supported by GA\v CR grant no.~21-27941S
 and RVO funding for I\v CO~67985840.}
\begin{abstract} We~describe the Dirichlet space of $M$-harmonic functions,
i.e.~functions annihilated by the invariant Laplacian on~the unit ball of the complex $n$-space,
as~the limit of the analytic continuation (in~the spirit of Rossi and Vergne) of the corresponding
weighted Bergman spaces. Characterizations in terms of tangential derivatives are given,
and the associated inner product is shown to be Moebius invariant.
The pluriharmonic and harmonic cases are also briefly treated.
\end{abstract}

\maketitle

\section{Introduction}
Let $\Bn$ be the unit ball in the complex $n$-space~$\CC^n$, $n\ge1$, and consider the standard
weighted Bergman spaces
$$ \AA_s(\Bn) := \{ f\in L^2(\Bn,d\mu_s): \;f\text{ is holomorphic on }\Bn\}  $$
of holomorphic functions on $\Bn$ square-integrable with respect to the measure
\[ d\mu_s(z) := \frac{\Gamma(s+n+1)}{\pi^n\Gamma(s+1)} (1-|z|^2)^s \, dz, \qquad s>-1,  \label{TA} \]
where $dz$ denotes the Lebesgue volume on~$\CC^n$. The~restriction on $s$ ensures that these
spaces are nontrivial, and the factor $\frac{\Gamma(s+n+1)}{\pi^n\Gamma(s+1)}$ makes $d\mu_s$
a~probability measure, so~that the function $\jedna$ (constant~one) has unit norm. It~is well known
that $\AA_s$ is a reproducing kernel Hilbert space, with reproducing kernel
\[ K^\hol_s(x,y) = (1-\spr{x,y})^{-s-n-1} . \label{TB}  \]
In~terms of the Taylor coefficients, a~holomorphic function $f(z)=\sum_\nu f_\nu z^\nu$ on $\Bn$
belongs to $\AA_s$ if and only~if
\[ \|f\|^2_s := \sum_\nu |f_\nu|^2 \frac{\nu!\Gamma(n+s+1)}{\Gamma(n+s+|\nu|+1)} < +\infty \label{TC}  \]
and the sum then coincides with the squared norm in~$\AA_s$. (Here the sum runs over all multi-indices~$\nu$,
and we are employing the usual multi-index notations.) Alternatively, in~terms of the homogeneous
components $f_m(z):=\sum_{|\nu|=m} f_\nu z^\nu$,
\[ \|f\|^2_s = \sum_{m=0}^\infty \frac{\Gamma(n+s+1)\Gamma(n+m)}{\Gamma(n)\Gamma(n+s+m+1)} \|f_m\|^2_\pBn
 = \sum_{m=0}^\infty \frac{(n)_m}{(n+s+1)_m} \|f_m\|^2_\pBn , \label{TD}  \]
where $\|f\|_\pBn$ denotes the norm in the space $L^2(\pBn,d\sigma)$ with respect to the normalized
surface measure $d\sigma$ on~$\pBn$. Here $(x)_m:=x(x+1)\dots(x+m-1)$ denotes the usual
Pochhammer symbol (rising factorial).

It~is a remarkable fact --- which prevails in the much more general context of bounded
symmetric domains, constituting the ``analytic continuation'' of the principal series
representations of certain semisimple Lie groups, cf.~Rossi and Vergne~\cite{RV} ---
that the weighted Bergman kernels $K^\hol_s(x,y)$, $s>-1$, continue to be positive definite
kernels in the sense of Aronszajn~\cite{Aro} for all $s\ge-n-1$, yielding thus an
``analytic continuation'' of the spaces~$\AA_s$. (One~calls the interval $\itv{-n-1,+\infty}$
the \emph{Wallach set} of~$\Bn$.) For $s>-n-1$, the norm in $\AA_s$ is still given by
\eqref{TC} and~\eqref{TD}. For $s=-n-1$, the kernel \eqref{TB} becomes constant~one,
and the corresponding reproducing kernel Hilbert space thus reduces just to the constants.
However, a~much more interesting space arises as the ``residue'' of $\AA_s$ at $s=-n-1$:
namely, the~limit
\[ \lim_{s\searrow-n-1} \frac{K^\hol_s(x,y)-1}{s+n+1} = \log\frac1{1-\spr{x,y}} =: K^\hol_\circ(x,y) \label{TE} \]
is~a~positive definite kernel on $\Bn\times\Bn$, and the associated reproducing kernel
Hilbert space --- denoted $\AA_\circ$ --- consists of all $f$ holomorphic on~$\Bn$ for which
\[ \begin{aligned}
\|f\|^2_\circ &:= \sum_m \lim_{s\searrow-n-1} \frac{(n)_m}{(n+s+1)_m} \|f_m\|^2_\pBn \\
 &= \sum_m m\frac{(n)_m}{m!} \|f_m\|^2_\pBn = \sum_\nu |\nu|\frac{\nu!}{|\nu|!} |f_\nu|^2 < +\infty,
\end{aligned}  \label{TF} \]
and $\|f\|_\circ$ gives the semi-norm on~$\AA_\circ$. This space is nothing else but the
familiar \emph{Dirichlet space} on~$\Bn$, see e.g.~Chapter~6.4 in~Zhu~\cite{Zhu},
where it is furthermore shown that the space $\AA_\circ$ and the above semi-inner product
are \emph{Moebius invariant}, in~the sense that $f\in\AA_\circ\implies f\circ\phi\in\AA_\circ$ and
$$ \spr{f,g}_\circ = \spr{f\circ\phi,g\circ\phi}_\circ  $$
for any biholomorphic self-map $\phi$ of~$\Bn$.

The~goal of the present paper is to exhibit the \Mh analogue of the construction above.

Recall that a function on $\Bn$ is~called \emph{Moebius-harmonic} (or~\emph{invariantly harmonic}),
or \emph{\Mh }for short, if~it is annihilated by the invariant Laplacian
\[ \td = 4(1-|z|^2) \sum_{j,k=1}^n (\delta_{jk}-z_j\oz_k)\frac{\partial^2}{\partial z_j\partial\oz_k}.
 \label{TH} \]
It~is standard (see e.g.~Rudin~\cite{Ru}, Stoll~\cite{Stoll}, or Chapter~6 in Krantz~\cite{KraPDECA})
that $\td$ commutes with biholomorphic self-maps (Moebius maps) of the ball:
$$ \td(f\circ\phi) = (\td f)\circ\phi, \qquad \forall f\in C^2(\Bn), \phi\in\autBn; $$
and, accordingly, $f$~is \Mh if and only if $f\circ\phi$~is, for any $\phi\in\autBn$.
The~class of \Mh functions lies in a way on the crossroads between the holomorphic and the harmonic functions:
it~resembles the latter in the sense that it is preserved by complex conjugation, while resembling the
former by reflecting the complex structure inherent in the invariance of the Laplacian~$\td$.
One~again has the \emph{\Mh weighted Bergman spaces}
$$ \MM_s := \{f\in L^2(\Bn,d\mu_s):\;f\text{ is \Mh on }\Bn\}   $$
which are nontrivial if and only if $s>-1$. The~role of the Taylor coefficients or, rather,
homogeneous components from the holomorphic case is now played by the decomposition into
``bi-graded spherical harmonics''. Namely, under the action of the group $U(n)$ of unitary
linear maps of~$\CC^n$, the~space $L^2(\pBn,d\sigma)$ decomposes into irreducible components
\[ L^2(\pBn,d\sigma) = \bigoplus_{p,q=0}^\infty \chpq,  \label{TI}  \]
where $\chpq$ is the space of restrictions to the sphere of harmonic polynomials on $\CC^n$
homogeneous of degree $p$ in $z$ and of degree $q$ in~$\oz$. Performing such a decompisition
on each sphere $|z|\equiv\text{const.}$ leads to the analogous Peter-Weyl decomposition
$$ \MM_s = \bigoplus_{p,q} \bhpq,  $$
where the space $\bhpq$ of ``solid harmonics''
$$ \bhpq = \{f\in C(\overline{\Bn}): \;f\text{ is \Mh on }\Bn \text{ and }f|_\pBn\in\chpq\} ,  $$
and the norm of $f=\sum_{p,q}f_{pq}$, $f_{pq}\in\bhpq$, is~given~by
\[ \|f\|^2_s = \sum_{p,q=0}^\infty C_{pq}(s) \|f_{pq}\|^2_\pBn,  \label{TJ}  \]
with the coefficients $C_{pq}(s)$ --- the counterparts of the $\frac{(n)_m}{(n+s+1)_m}$ from the
holomorphic case --- given by an explicit formula involving hypergeometric functions;
see Section~\ref{Sec2} below for the details. Finally, the~reproducing kernel of the space~$\MM_s$,
$s>-1$, is given~by
\[ K_s(z,w) = \sum_{p,q} \frac{K_{pq}(z,w)}{C_{pq}(s)} ,  \label{TK}  \]
where $K_{pq}(zw)$ is the reproducing kernel of~$\bhpq$ (with the inner product inherited from~$L^2(\pBn,d\sigma)$),
for~which there is again an explicit formula.
Now~it has been shown in Section~6.3 of~\cite{EY} that, exactly as in the holomorphic case,
$K_s(z,w)$ extends to a holomorphic function of $s$ on $\Re s>-n-1$, continues to be a positive
definite kernel on $\Bn\times\Bn$ for all $s\ge-n-1$ (and~only for these~$s$), and the norm in
the corresponding reproducing kernel Hilbert space --- the ``analytic continuation'' of~$\MM_s$ ---
for $s>-n-1$ is still given~by~\eqref{TJ}. (Thus the ``\Mh Wallach set' of $\Bn$ is again the
interval $\itv{-n-1,+\infty}$.)

Motivated by the considerations for the holomorphic case, we~are now interested in the limit
as $s\searrow-n-1$. In~contrast to the holomorphic case, this can now be done in \emph{three} ways.

\begin{itemize}
\item[(a)] We~simply take $s=-n-1$ in~\eqref{TK}. The kernel $K_s(z,w)$ reduces to constant~one,
and the corresponding reproducing kernel Hilbert space thus again reduces just to the constants,
with $\|\jedna\|=1$. This is the trivial case.
\item[(b)] As~in the holomorphic case, next we take
\[ \lim_{s\searrow-n-1} \frac{K_s(z,w)-1}{n+s+1} = \log\frac1{|1-\spr{x,y}|^2} =: K_\circ(x,y). \label{TL} \]
This is a positive definite kernel, and the corresponding reproducing kernel Hilbert space ---
denoted $\MM_\circ$ --- consists precisely of the orthogonal sum of the holomorphic Dirichlet space
$\AA_\circ$ above and its complex conjugate $\overline{\AA_\circ}$ (the~constants being counted,
of~course, only~once), with the (semi-)norm given~by
\[ \|f+\overline g\|^2_\circ := \|f\|_\circ^2 +\|g\|^2_\circ.  \label{TM}  \]
In~some sense, one~can perhaps view $\MM_\circ$ as the \emph{\Ph Dirichlet space}.
\item[(c)] Finally, we~can take
\[ \lim_{s\searrow-n-1} \frac{K_s(x,y)-1-(n+s+1)K_\circ(x,y)}{(n+s+1)^2} =: K_\cici(x,y) , \label{TN}  \]
which is a positive-definite kernel on $\Bn\times\Bn$, with the (semi-)norm in the corresponding
reproducing kernel Hilbert space --- denoted $\MM_\cici$ --- given~by
\[ \|f\|^2_\cici := \sum_{p,q} \lim_{s\searrow-n-1} (n+s+1)^2 C_{pq}(s) \|f_{pq}\|^2_\pBn \label{TO}  \]
(hence, this time, all \Ph functions get zero norm). This~is, by~definition,
the \emph{\Mh Dirichlet space}.
\end{itemize}

The~occurrence of case (c) arises from the fact that the coefficient functions $C_{pq}(s)$ now
turn out to have a double pole at $s=-n-1$ (in~contrast to the single pole of $\frac{(n)_m}{(n+s+1)_m}$
in the holomorphic situation). Again, this phenomenon is clearly reminiscent of the ``composition series''
arising in the theory of analytic continuation of holomorphic discrete series representations
mentioned above, cf.~Faraut and Koranyi~\cite{FKjfa}.
Note also that the case (c) disappears completely when $n=1$; thus $\MM_\cici$ is relevant only for
$\Bn$ with $n\ge2$.

Our~first result is a direct formula for the semi-norm in~$\MM_\cici$.

\begin{corollary*} {\rm(Corollary~\ref{PX})}
In~terms of the Peter-Weyl decomposition $f=\sum_{p,q}f_{pq}$, $f_{pq}\in\bhpq$,
$$ \|f\|^2_\cici = \sum_{p,q} \frac{(p)_n(q)_n}{\Gamma(n)^2} \|f_{pq}\|^2_{\pBn}.  $$
\end{corollary*}

By~polarization, this of course implies also the corresponding formula
\[ \spr{f,g}_\cici = \sum_{p,q} \frac{(p)_n(q)_n}{\Gamma(n)^2} \spr{f_{pq},g_{pq}}_{\pBn}  \label{TP}  \]
for the semi-inner product in~$\MM_\cici$.

Next we give a description of $\MM_\cici$ that does not involve the Peter-Weyl components.
It~is easy to see that the averaging operator
$$ \Pi_0 f(r\zeta) := \int_\pBn f(r\eta) \, d\sigma(\eta), \qquad 0\le r<1, \zeta\in\pBn,  $$
(this can also be written as $\Pi_0 f(z)=\int_{U(n)} f(kz)\,dk$, where $dk$ stands for the normalized
Haar measure on the compact group~$U(n)$) is~just the projection
$$ f\longmapsto f(0)\jedna $$
of \Mh functions onto the subspace $\mathbf H^{00}$ of constants. Similarly, one~can give projections
$\Pi$ and $\overline\Pi$ onto the subspaces $\bigoplus_p \mathbf H^{p0}$ and $\bigoplus_q\mathbf H^{0q}$
of~the holomorphic and anti-holomorphic functions, respectively (explicit formulas for $\Pi$ and $\overline\Pi$
will be given in Section~\ref{Sec2} below). Hence, we~also have
$$ P := \Pi+\overline\Pi-\Pi_0,  $$
the projection onto the subspace of \Ph functions, and
$$ Q := I-P,  $$
the projection onto their orthogonal complement $\bigoplus_{p,q\ge1}\bhpq$.

(All these projections are automatically orthogonal with respect to any $U(n)$-invariant inner product,
but make sense also in complete generality on the vector space $\MM$ of all \Mh functions on~$\Bn$.)

Consider now the tangential vector fields
$$ L_{jk} := \oz_j \partial_k - \oz_k \partial_j, \qquad
 \overline L_{jk} := z_j \dbar_k - z_k \dbar_j, \qquad 1\le j,k\le n, \quad j\neq k,  $$
and denote by $\LL_m$, $1\le m\le2n(n-1)$, the collection of all these operators
(in~some fixed order). Finally, for any function $f$ on~$\Bn$ and $0<r<1$, let $f_r$ be the
function on $\pBn$ defined by
$$ f_r(\zeta) := f(r\zeta) $$
and denote
$$ \|f\|^2_\Hardy := \sup_{0<r<1} \|f_r\|^2_\pBn.  $$

\begin{theorem*} {\rm(Theorem~\ref{PF})}
If~$f$ is \Mh on $\Bn$, $n\ge2$, then $f\in\MM_\cici$ if and only if
$$ \sum_{j_1,j_2,\dots,j_n=1}^{2n(n-1)} \|\LL_{j_1}\LL_{j_2}\dots\LL_{j_n}(Qf)\|^2_\Hardy <+\infty, $$
and the square root of the left-hand side is a seminorm equivalent to~$\|f\|_\cici$.
\end{theorem*}

As~a~consequence, we~obtain the following complete analogues of the holomorphic case.

\begin{corollary*} {\rm(Corollary~\ref{PG})}
The space $\MM_\cici$ is Moebius invariant: $f\in\MM_\cici$ implies $f\circ\phi\in\MM_\cici$
for any $\phi\in\autBn$.
\end{corollary*}

\begin{theorem*} {\rm(Theorem~\ref{PK})}
The~semi-inner product \eqref{TP} is Moebius invariant:
$$ \spr{f,g}_\cici = \spr{f\circ\phi,g\circ\phi}_\cici \qquad \forall f,g\in\MM_\cici, \forall\phi\in\autBn.  $$
\end{theorem*}

The~usual proof of the analogue of the last corollary for the holomorphic case (cf.~Zhu~\cite{Zhu}, Theorem~6.13)
relies on the use of radial derivatives and their generalizations; this approach unfortunately breaks down in
the \Mh situation. Similarly, the usual proof of the last theorem in the holomorphic case relies on explicit
computations involving Taylor coefficients (cf.~Zhu~\cite{Zhu}, Theorem~6.15), which becomes hopeless for~$\MM_\cici$;
we~instead use an argument employing analytic continuation. In~both cases, our~approach here can be used to give
a new description of the classical Dirichlet space on the ball and a new proof of the invariance of the Dirichlet
inner product in the holomorphic case.

Additionally, the~methods just mentioned apply also to the \Ph Dirichlet space~$\MM_\circ$,
which seems to have received basically no attention at all in the literature.
The~second result below apparently has no counterpart in the \Mh case.

\begin{theorem*} {\rm(Theorem~\ref{PH})}
If $f$ is \Ph on $\Bn$, $n\ge2$, then $f\in\MM_\circ$ if and only if
$$ \sum_{j_1,j_2,\dots,j_n=1}^{2n(n-1)} \|\LL_{j_1}\LL_{j_2}\dots\LL_{j_n}f\|^2_\Hardy <+\infty $$
if and only if
$$ \sum_{j_1,j_2,\dots,j_{n+k+1}=1}^{2n(n-1)} \|\LL_{j_1}\LL_{j_2}\dots\LL_{j_{n+k+1}}f\|^2_k <+\infty, $$
for some (equivalently, any) nonnegative integer~$k$.

Furthermore, the square root of the above quantities is a seminorm equivalent to~$\|f\|_\circ$.
\end{theorem*}

Here $\|\cdot\|_s$ denotes, more generally, the~norm in $L^2(\Bn,d\mu_s)$, $s>-1$ (i.e.~not only on its
holomorphic or \Mh subspaces).

Let
$$ \cN := \sum_{j=1}^n z_j\partial_j + \oz_j\dbar_j  $$
denote the radial derivative operator on~$\Bn$; note than $\cN f$ is \Ph whenever $f$~is.

\begin{theorem*} {\rm(Theorem~\ref{PI})}
If $f$ is \Ph on $\Bn$, $n\ge1$, then $f\in\MM_\circ$ if and only if
$$ \|\cN^m f\|^2_{2m-n-1} <+\infty  $$
for some (equivalently, any) integer~$m>\frac n2$.
Furthermore, the square root of the left-hand side is a seminorm equivalent to~$\|f\|_\circ$.
\end{theorem*}

For $f$ holomorphic the last two theorems, of~course, give criteria for $f$ to belong to the
ordinary Dirichlet space $\AA_\circ$ on~$\Bn$; the~second one is then common knowledge,
and though the first of them must surely also be known to experts in the field,
the~authors were unable to pinpoint a specific reference to the literature.

Finally, for~the sake of completeness, we~give details also for the context of \emph{harmonic} functions,
where the corresponding \emph{harmonic Dirichlet space} appears in the literature under various names.
Namely, consider this time the unit ball $\bn$ of~$\RR^n$, $n\ge2$, and for any $s>-1$ let
$$ \HH_s(\bn) := \{ f\in L^2(\bn,d\rho_s): \;f\text{ is harmonic on }\bn\}  $$
be~the \emph{weighted harmonic Bergman space} of all harmonic functions on $\bn$ square-integrable
with respect to the measure
$$ d\rho_s(x) := \frac{\Gamma(\frac n2+s+1)}{\pi^{n/2}\Gamma(s+1)} (1-|x|^2)^s\,dx,  $$
where $dx$ denotes the Lebesgue volume on~$\RR^n$. The~restriction on $s$ ensures that these spaces are nontrivial,
and the factor $\frac{\Gamma(\frac n2+s+1)}{\pi^{n/2}\Gamma(s+1)}$ makes $d\rho_s$ a probability measure,
so~that $\|\jedna\|=1$. For~each $p\ge0$, denote by $\bhp$ the space of harmonic polynomials on $\RR^n$
homogeneous of degree~$p$. Any~harmonic function $f$ on $\bn$ then admits a (unique) decomposition
\[ f = \sum_{p=0}^\infty f_p, \qquad f_p\in\bhp,  \label{TV}  \]
and the space $\HH_s$ decomposes~as
\[ \HH_s = \bigoplus_{p=0}^\infty \bhp , \label{TW} \]
with the norm given~by
\[ \|f\|^2_s = \sum_{p=0}^\infty \frac{(\frac n2)_p}{(\frac n2+s+1)_p} \|f_p\|_\pbn , \label{TX} \]
where $\|\cdot\|_\pbn$ denotes the norm in $L^2(\pbn,d\sigma)$ with respect to the normalized
surface measure $d\sigma$ on the unit sphere~$\pbn$. The~reproducing kernel of~$\HH_s$ is given~by
\[ K^\harm_s(x,y) = \sum_p \frac{(\frac n2+s+1)_p}{(\frac n2)_p} Z_p(x,y),  \label{TY}  \]
where the \emph{zonal harmonic} $Z_p(x,y)$ is the reproducing kernel of~$\bhp$
(with respect to the norm $\|\cdot\|_\pbn$). It~follows that $K^\harm_s(x,y)$ extends as a holomorphic
function of $s$ to the entire complex plane, and continues to be a positive definite kernel on $\bn\times\bn$
for any $s\ge-\frac n2-1$ (and only for these~$s$). For $s>-\frac n2-1$, the~norm in the corresponding
reproducing kernel Hilbert spaces (still denoted by~$\HH_s$) is still given by the formula~\eqref{TX}.
For $s=-\frac n2-1$, \eqref{TY} again reduces just to constant~one, and the corresponding reproducing
kernel Hilbert space thus reduces to the constants; while the limit
$$ \lim_{s\searrow-\frac n2-1} \frac{K^\harm_s(x,y)-1}{s+\frac n2+1} =: K^\harm_\sq(x,y)  $$
is a positive definite kernel on~$\bn\times\bn$, corresponding to the reproducing kernel Hilbert space
with the (semi)norm given~by
\[ \|f\|^2_\sq := \sum_p \lim_{s\searrow-\frac n2-1} (s+\tfrac n2+1)\frac{(\frac n2)_p}{(\frac n2+s+1)_p} \|f_p\|_\pbn
 = \sum_p p\frac{(\frac n2)_p}{p!} \|f_p\|_\pbn .  \label{TZ}  \]
This space --- denoted $\HH_\sq$ --- is~the \emph{harmonic Dirichlet space} on~$\bn$.

Let $X_{jk}$, $j,k=1,\dots,n$, $j\neq k$, denote the tangential vector fields
$$ X_{jk} = x_j \partial_k - x_k \partial_j  $$
on~$\RR^n$, and denote by $\XX_m$, $1\le m\le n(n-1)$, the collection of all these operators
(in~some fixed order).

\begin{theorem*} {\rm(Theorem~\ref{PJ})}
If $f$ is harmonic on~$\bn$, $n\ge2$, then $f\in\HH_\sq$ if and only~if
$$ \sum_{j_1,\dots,j_m=1}^{n(n-1)} \|\XX_{j_1}\dots\XX_{j_m}f\|^2_{2m-\frac n2-1} <+\infty  $$
for some (equivalently, any) integer $m>\frac n4$.
Furthermore, the square root of the left-hand side is a seminorm equivalent to~$\|f\|_\sq$.
\end{theorem*}

This time there is no Moebius invariance, since the Moebius self-maps of $\bn$ do not preserve
harmonicity for $n>2$.

We~review the necessary background material in Section~\ref{Sec2}, then present the basic properties
of the \Mh Dirichlet space $\MM_\cici$ in Section~\ref{Sec3}. Moebius invariance is discussed in Section~\ref{Sec4}.
The~\Ph Dirichlet space is treated in Section~\ref{Sec5}, and the harmonic case in Section~\ref{Sec6}.

Throughout the paper, the notation
$$ A \asymp B  $$
means that
$$ cA \le B \le \frac1c A  $$
for some $0<c\le1$ independent of the variables in question. The~symbols $\frac\partial{\partial z_j}$
and $\frac\partial{\partial\oz_j}$, commonly abbreviated just to $\partial_j$ and $\dbar_j$, respectively,
stand for the usual Wirtinger operators on~$\CC^n$; similarly on~$\RR^n$,
$\partial_k$ stands for $\frac\partial{\partial x_k}$.
For~typesetting reasons, the inner product $\spr{x,y}$ in $\CC^n$ is sometimes also denoted by~$\cdd xy$.
Finally, $\ZZ,\NN,\RR$ and $\CC$ denote the sets of all integers, all nonnegative integers, all real
and all complex numbers, respectively.

\section{Notation and preliminaries} \label{Sec2}
The stabilizer of the origin $0\in\Bn$ in $\autBn$ is the group $\Un$
of all unitary transformations of~$\CC^n$; that~is, of~all linear operators $U$ that
preserve inner products:
$$ \spr{Uz,Uw} = \spr{z,w} \qquad \forall z,w\in\CC^n.   $$
Each $U\in\Un$ maps the unit sphere $\pBn$ onto itself, and the surface measure
$d\sigma$ on $\pBn$ is invariant under~$U$. It~follows that the composition with
elements of~$\Un$,
\[ T_U: f \mapsto f \circ U^{-1},  \label{tVA}  \]
is~a unitary representation of $\Un$ on $L^2(\pBn,d\sigma)$. We~will need the decomposition
of this representation into irreducible subspaces. These turn out to be given by
\emph{bigraded spherical harmonics}~$\chpq$; the standard sources for this are
Rudin~\cite[Sections 12.1--12.2]{Ru}, or Krantz~\cite[Sections 6.6--6.8]{KraPDECA},
with basic ingredients going back to Folland~\cite{Foll}.

Namely, for integers $p,q\ge0$, let $\chpq$ be vector space of restrictions to $\pBn$
of harmonic polynomials $f(z,\oz)$ on $\CC^n$ which are homogeneous of degree $p$ in $z$
and homogeneous of degree $q$ in~$\oz$. Then $\chpq$ is invariant under the action
\eqref{tVA} of~$\Un$, is $\Un$-irreducible (i.e.~has no proper $\Un$-invariant subspace)
and
\[ L^2(\pBn,d\sigma) = \bigoplus_{p,q=0}^\infty \chpq .  \label{tVB}  \]
Furthermore, if $T$ is a linear operator on $L^2(\pBn,d\sigma)$ commuting with the action~\eqref{tVA}
i.e.~$T(f\circ U)=(Tf)\circ U$ for all $U\in\Un$ and $f\in L^2(\pBn,d\sigma)$),
then $T$ is diagonalized by the decomposition~\eqref{tVA}, i.e.~$T$ maps each $\chpq$
into itself and $T|\chpq=c_{pq}I|\chpq$ for some complex constants $c_{pq}$,
where $I$ denotes the identity operator.

Since each space $\chpq$ is finite-dimensional, the evaluation functional $f\mapsto f(\zeta)$
at each $\zeta\in\pBn$ is automatically continuous on~it; it~follows that $\chpq$ --- with the
inner product inherited from $L^2(\pBn,d\sigma)$ --- has a reproducing kernel. This reproducing
kernel turns out to be given by $H^{pq}(\cdd\zeta\eta)$, where for $n\ge2$
\[ \begin{aligned}
H^{pq}(z) &= \frac{(-1)^q(n+p+q-1)(n+p-2)!}{(n-1)!q!(p-q)!}  \\
&\hskip4em \,\times z^{p-q} \FF21{-q,n+p-1}{p-q+1}{|z|^2} \qquad\text{for }p\ge q, \end{aligned} \label{tVE}  \]
while $H^{pq}(z)=H^{qp}(\oz)$ for $p<q$.
For $n=1$, the spaces $\chpq$ reduce just to $\{0\}$
if $pq\neq0$, while $\mathcal H^{p0}=\CC z^p$, $\mathcal H^{0q}=\CC\oz^q$ and
$H^{p0}(z)=z^p$, $H^{0q}(z)=\oz^q$; note that the formula \eqref{tVE} still works for $n=1$ and $pq=0$.

Denote
\[ \begin{aligned}
S^{pq}(r) &:= r^{p+q} \FF21{p,q}{p+q+n}{r^2} \Big/ \FF21{p,q}{p+q+n}1 \\
 &= \frac{\Gamma(p+n)\Gamma(q+n)}{\Gamma(n)\Gamma(p+q+n)} r^{p+q} \FF21{p,q}{p+q+n}{r^2} . \end{aligned}
 \label{tVF}  \]
Then for each $f\in\chpq$, the (unique) solution to the Dirichlet problem $\td u=0$ on~$\Bn$,
$u|_\pBn=f$ is given~by
\[ u(r\zeta) = S^{pq}(r) f(\zeta), \qquad 0\le r\le 1, \; \zeta\in\pBn.  \label{tVG}  \]
Many of the formulas above originate in~\cite{Foll}.

For each $p,q\ge0$, let $\bhpq$ be the space of all functions on $\Bn$ of the form
\eqref{tVG} with $f\in\chpq$. In~other words, while $\chpq$ is the space of spherical
harmonics on the sphere~$\pBn$, $\bhpq$ is the associated space of ``solid'' \Mh functions
on~$\Bn$. With the inner product inherited from~$L^2(\pBn,d\sigma)$, each $\bhpq$ is thus
a finite-dimensional Hilbert space of \Mh functions on~$\Bn$, unitarily isomorphic to
the space $\chpq$ via the isomorphism~\eqref{tVG}, and with reproducing kernel
\[ K^{pq}(r\zeta,R\xi) := S^{pq}(r) S^{pq}(R) H^{pq}(\cdd\zeta\xi). \label{tVM}  \]

It~was shown in Proposition~1 in~\cite{EY} that if $H$ is any Hilbert space of \Mh functions
on $\Bn$ whose inner product is invariant under rotations (i.e.~$f\in H$, $U\in\Un$ imply
$f\circ U\in H$ and $\|f\circ U\|_H=\|f\|_H$), then the following is true:
\begin{itemize}
\item[--] for each~$p,q\ge0$, $H\cap\bhpq$ is either $\{0\}$ or all of~$\bhpq$;
\item[--] if $\bhpq\subset H$, then
$$ \spr{f,g}_H = c_{pq} \spr{f,g}_\pBn \qquad \forall f,g\in\bhpq $$
with some constant $c_{pq}$, $0<c_{pq}<+\infty$ (independent of~$f,g$);
\item[--] the reproducing kernel of $H$ is given by
\[ K_H(x,y) = \sum_{p,q:\bhpq\subset H} \frac{K_{pq}(x,y)}{c_{pq}},  \label{tA}  \]
with the sum converging locally uniformly on $\Bn\times\Bn$,
as~well as in~$H$ as a function of $x$ for each fixed $y$, or~vice versa.
\end{itemize}

One~can also formally define $c_{pq}:=+\infty$ if $\bhpq\not\subset H$; then $H$ contains precisely
those $\bhpq$ for which $c_{pq}$ is finite, and in \eqref{tA} the summation can be extended over all~$p,q\ge0$,
with the usual convention that $1/\infty:=0$.

One~can also allow semi-Hilbert spaces, i.e.~with semi-definite (semi-)norm instead of norm;
then the above still holds, except that $c_{pq}$ can be zero for some $p,q$ and \eqref{tA} is the
reproducing kernel not for $H$ but only for the subspace $H_0:=\bigoplus\{\bhpq: c_{pq}>0\}$.

For the weighted \Mh Bergman spaces of all \Mh functions in $L^2(\Bn,(1-|z|^2)^s\,dz)$, $s>-1$,
an~explicit formula for the $c_{pq}=: c_{pq}(s)$ was given in~(69) in~\cite{EY}. Renormalizing so
as to pass to our normalized measures from~\eqref{TA}, and our spaces $\MM_s$ from the Introduction,
it~was shown in Section~6.3 of~\cite{EY} that the corresponding constants $C_{pq}(s)=c_{pq}(s)/c_{00}(s)$
are given~by
\[ C_{pq}(s) = \frac1{\Gamma(n)} \int_0^1 G^{(n)}_{pq}(t) \; (1-t)^{n+s} \,dt , \label{tC}  \]
where $G_{pq}$ is the function
$$ G_{pq}(t) := \frac{\Gamma(n+p)^2\Gamma(n+q)^2}{\Gamma(n)^2\Gamma(n+p+q)^2} t^{p+q+n-1}
 \FF21{p,q}{n+p+q}t ^2  $$
(see formula (98) there). Furthermore, it~was shown there that $G_{pq}^{(n)}$ is positive and continuous
on $(0,1)$, with a finite value at the origin and $G^{(n)}_{pq}(t)=O(\log\frac1{1-t})$ as $t\nearrow1$,
from which it follows that \eqref{tC} furnishes an analytic continuation of $C_{pq}(s)$ to $\Re s>-n-1$
and $0<C_{pq}(s)<+\infty$ for $s\in(-n-1,+\infty)$. The~functions
$$ K_s(z,w) = \sum_{p,q} \frac{K_{pq}(z,w)}{C_{pq}(s)}  $$
thus continue to be positive definite kernels on $\Bn\times\Bn$ in the sense of Aronszajn~\cite{Aro},
and the norm in the associated reproducing kernel Hilbert spaces --- still denoted~$\MM_s$ ---
is~given~by
$$ \|f\|^2_s = \sum_{p,q} C_{pq}(s) \|f_{pq}\|^2_\pbn, \qquad s>-n-1,  $$
for an \Mh function $f$ with Peter-Weyl decomposition $f=\sum_{p,q}f_{pq}$, $f_{pq}\in\bhpq$.

We~conclude this section by describing the projections onto the Peter-Weyl components~$\bhpq$.
From the reproducing property of $H^{qp}(\cdd\zeta\eta)$ and~\eqref{tVG}, we~have,
for any \Mh function $f$ on~$\Bn$,
\[ f_{pq}(r\zeta) = \int_\pbn f(r\eta) H^{pq}(\cdd\zeta\eta) \,d\sigma(\eta)
 \qquad \forall\zeta\in\pbn, \forall r\in(0,1).  \label{UX}  \]
In~particular, for $p=q=0$,
\[ \Pi_0 f(r\zeta) = \int_\pbn f(r\eta) \,d\sigma(\eta), \qquad 0\le r<1, \;\zeta\in\pBn, \label{UZ} \]
is the projection onto the constants~$\mathbf H^{00}$.

\begin{proposition} \label{PA}
For any $f$ \Mh on~$\Bn$, the limit
\[ \Pi f(r\zeta) := \lim_{R\nearrow1} \int_\pbn \frac{f(r\eta)}{(1-R\spr{\zeta,\eta})^n} \,d\sigma(\eta) \label{UY} \]
exists, and $\Pi f$ equals the projection of $f$ onto the subspace $\AA:=\bigoplus_{P=0}^\infty\mathbf H^{p0}$
of holomorphic functions on~$\Bn$.
\end{proposition}

\begin{proof} Expanding $(1-R\spr{\zeta,\eta})^{-n}$ by the binomial formula shows that the integral equals
\begin{align*}
& \int_\pbn f(r\eta) \sum_{j=0}^\infty R^j \frac{(n)_j}{j!} \spr{\zeta,\eta}^j \, d\sigma(\eta) \\
& \hskip4em = \int_\pbn f(r\eta) \sum_{j=0}^\infty R^j H^{j0}(\cdd\zeta\eta) \, d\sigma(\eta) \qquad\text{by \eqref{tVE}} \\
& \hskip4em = \sum_j R^j f_{j0}(r\zeta) \qquad\text{by \eqref{UX}}  \\
& \hskip4em = \sum_j f_{j0}(Rr\zeta) = \Pi f(Rr\zeta),
\end{align*}
the interchange of the summation and integration signs being justified by the locally uniform convergence.
Letting $R\nearrow1$, the claim follows.
\end{proof}

We~remark that it~is, of~course, not~possible to interchange the limit and the integral in~\eqref{UY},
since $(1-\spr{\zeta,\cdot})^{-n}$ is not integrable over~$\pBn$, for any $\zeta\in\pBn$.

Taking complex conjugates in~\eqref{UY}, one gets also the projection $\overline\Pi$ onto the subspace
of anti-holomorphic functions, and the projection
$$ P := \Pi + \overline\Pi - \Pi_0  $$
onto the subspace $\cP:=\bigoplus_{pq=0}\bhpq$ of \Ph functions on~$\Bn$.

For~later use, we denote~by
$$ Q := I-P  $$
the projection onto the orthogonal complement $\bigoplus_{pq>0}\bhpq$
(``the \Mh functions with no \Ph component'').

\section{\Mh Dirichlet space}  \label{Sec3}
For an \Mh function $f=\sum_{p,q}f_{pq}$, $f_{pq}\in\bhpq$, on~$\Bn$, we~have by \eqref{tVA} and~\eqref{tVG}
\[ \|f\|^2_s = \sum_{p,q} C_{pq}(s) \|f_{pq}\|^2_\pBn  \label{VA}  \]
for any $s>-1$, with
\begin{align}
C_{pq}(s) &:= \frac{\Gamma(n+s+1)}{\pi^n\Gamma(s+1)} \int_0^1 S^{pq}(r)^2 \frac{2\pi^n}{\Gamma(n)} (1-r^2)^s r^{2n-1}\,dr \nonumber \\
&= \frac{(s+1)_n}{\Gamma(n)} \int_0^1 S^{pq}(\sqrt t)^2 t^{n-1}\,(1-t)^s\, dt \nonumber \\
&= \frac{(s+1)_n}{\Gamma(n)} \int_0^1 G_{pq}(t) \,(1-t)^s \,dt, \label{VB}
\end{align}
where $G_{pq}$ is the function
\[ G_{pq}(t) := t^{p+q+n-1} \frac{\Gamma(n+p)^2\Gamma(n+q)^2}{\Gamma(n)^2\Gamma(n+p+q)^2} \FF21{p,q}{n+p+q}t ^2. \label{VC}  \]
(See \cite{EY} for the details.)

The~content of the following lemma is standard; we~include the short proof for the sake of completeness.

\begin{lemma} \label{PB}
Let $F(z)=\sum_{k=0}^\infty F_kz^k$ be a holomorphic function on the disc $|z|<R$, $R>0$.
Then for any $\delta\in(0,R)$,
\begin{itemize}
\item[(a)] the integral
$$ \mathcal I(s) := \int_{1-\delta}^1 F(1-t)(1-t)^s \,dt, \qquad s>-1,  $$
extends to a holomorphic function of $s$ on the entire complex plane~$\CC$,
except for possible simple poles at $s=-j-1$, $j=0,1,2,\dots$, with residues~$F_j$;
\item[(b)] for $m=1,2,\dots$, the integral
$$ \mathcal I_m(s) := \int_{1-\delta}^1 F(1-t) \Big(\log\frac1{1-t}\Big)^m (1-t)^s \,dt, \qquad s>-1,  $$
extends to a holomorphic function of $s$ on the entire complex plane~$\CC$,
except for possible poles of multiplicity $m+1$ at $s=-j-1$, $j=0,1,2,\dots$,
of strength~$m!F_j$.
\end{itemize}
\end{lemma}

\begin{proof} (a) From the Taylor expansion, we~have for any $N=0,1,2,\dots$,
$$ F(z) = F_0 + F_1 z + \dots + F_{N-1} z^{N-1} + z^N G_N(z),  $$
with $G_N$ holomorphic on~$|z|<R$. By~uniform convergence,
\[ \mathcal I(s) = \sum_{j=0}^{N-1} F_j \frac{\delta^{s+j+1}}{s+j+1}
 + \int_{1-\delta}^1 G_N(1-t) (1-t)^{s+N} \, dt.  \label{VD}  \]
The $j$-th summand in the sum is holomorphic on $\CC$ except for a simple pole
at $s=-j-1$ with residue~$F_j$, while the integral is a holomorphic function
on $\Re s>-N-1$. As~$N$ was arbitrary, the~claim follows.

(b) Differentiating \eqref{VD} $m$ times with respect to $s$ yields
\[ \mathcal I_m(s) = \sum_{j=0}^{N-1} F_j \frac{m!\delta^{s+j+1}}{(s+j+1)^{m+1}}
 + \int_{1-\delta}^1 G_N(1-t) \Big(\log\frac1{1-t}\Big)^m (1-t)^{s+N} \, dt,  \label{VE}  \]
and the claim again follows.
\end{proof}

Note from \eqref{VB} that $C_{00}(s)\equiv1$ has no poles, while
\[ C_{p0}(s) = C_{0p}(s) = \frac{(n)_p}{(n+s+1)_p}  \label{VP}  \]
has simple poles at $s=-n-1,\dots,-n-p$.

\begin{proposition} \label{PC}
Let $pq>0$. Then $C_{pq}(s)$ extends to a holomorphic function of $s$ on the entire~$\CC$,
except for double poles at $s=-n-1,-n-2,\dots,-2n$ and triple poles at $s=-2n-1-j$,
$j=0,1,2,\dots$. The~double pole at $s=-n-1$ has strength $(p)_n(q)_n/\Gamma(n)^2$.
\end{proposition}

\begin{proof} Since $G_{pq}$ is continuous on the unit disc, the integral
$$ \int_0^{1-\delta} G_{pq}(t) (1-t)^s \,dt  $$
is a holomorphic function of $s$ on the entire complex plane, for any $0<\delta<1$.
For~the integral from $1-\delta$ to~$1$, formula (12) in \cite[\S2.10]{BE} tells us that
$$ \frac{\Gamma(n+p)\Gamma(n+q)}{\Gamma(n)\Gamma(n+p+q)} \FF21{p,q}{n+p+q}t
 = A_0(1-t) + A_1(1-t) (1-t)^n \log\frac1{1-t}  $$
for $|\arg(1-t)|<\pi$, with
\begin{align*}
A_0(z)&:= \sum_{j=0}^{n-1} \frac{(p)_j(q)_j}{(1-n)_j j!} z^j
 + \sum_{j=0}^\infty \frac{(-1)^n\Gamma(n+p+j)\Gamma(n+q+j)}{j!(j+n)!\Gamma(n)\Gamma(p)\Gamma(q)} z^{n+j}h''_j, \\
A_1(z) &:= \sum_{j=0}^\infty \frac{(-1)^n\Gamma(n+p+j)\Gamma(n+q+j)}{j!(j+n)!\Gamma(n)\Gamma(p)\Gamma(q)} z^j
\end{align*}
holomorphic on $|z|<1$; here
$$ h''_j := \psi(j+1)+\psi(j+n+1)-\psi(j+n+p)-\psi(j+n+q)  $$
where $\psi:=\Gamma'/\Gamma$ is the digamma function. It~follows that
$$ G_{pq}(t) = B_0(1-t) + B_1(1-t) (1-t)^n \log\frac1{1-t} + B_2(t) (1-t)^{2n} \Big(\log\frac1{1-t}\Big)^2 ,  $$
with
\begin{align*}
B_0(z) &= (1-z)^{p+q+n-1} A_0(z)^2, \\
B_1(z) &= 2(1-z)^{p+q+n-1} A_0(z) A_1(z), \qquad\text{and} \\
B_2(z) &= (1-z)^{p+q+n-1} A_1(z)^2
\end{align*}
holomorphic on $|z|<1$. Applying the lemma, we~see that
$$ \int_{1-\delta}^1 G_{pq}(t) (1-t)^s \, dt  $$
extends to an entire function of~$s$, except for possible simple poles at $s=-1,-2,\dots,-n$,
possible double poles at $s=-n-1,-n-2,\dots,-2n$ and possible triple poles at $s=-2n-1,-2n-2,\dots$.
The~strength of the double pole at $s=-n-1$ is
\[ 1!A_1(0) = \frac{(-1)^n\Gamma(n+p)\Gamma(n+q)}{n!\Gamma(n)\Gamma(p)\Gamma(q)}
 = \frac{(-1)^n}{n!\Gamma(n)} (p)_n (q)_n.  \label{VF}  \]
Finally, passing from $\int_0^1 G_{pq}(t)(1-t)^s\,dt$ to $C_{pq}(s)$, the~factor $\frac{(s+1)_n}{\Gamma(n)}$
cancels the simple poles at $s=-1,-2,\dots,-n$, while \eqref{VF} gets multiplied~by
$\frac{(s+1)_n}{\Gamma(n)}\Big|_{s=-n-1}=\frac{(-1)^nn!}{\Gamma(n)}$.
\end{proof}

\begin{remark} \label{XXCC}
A~variant of the last proof can be given by first integrating by parts to~get
\[ C_{pq}(s) = \frac1{\Gamma(n)} \int_0^1 G^{(n)}_{pq}(t) (1-t)^{n+s} \, dt  \label{XC}  \]
for any $s>-n-1$ and $p+q>0$; see Section~6.3 in~\cite{EY}. This has the advantage of showing
that $C_{pq}(s)$ is positive for any $s>-n-1$, and also is a decreasing function of~$s$
on~this interval. Since we are not interested in positivity or monotonicity at the moment,
it~was simpler to apply Lemma~\ref{PB} directly.  \qed  \end{remark}

Denoting the analytic continuation still by~$C_{pq}(s)$, we~thus see that
$$ C^\cici_{pq} := \lim_{s\searrow-n-1} (s+n+1)^2 C_{pq}(s)  $$
exists for all $p,q\ge0$, and equals
\[ C^\cici_{pq} = \begin{cases} 0 & \text{if }pq=0 \\
 \dfrac{(p)_n(q)_n}{\Gamma(n)^2} \qquad & \text{if }pq>0. \end{cases} \label{VG}  \]
We~thus arrive at the following corollary.

\begin{corollary} \label{PX}
In~terms of the Peter-Weyl decomposition $f=\sum_{p,q}f_{pq}$, $f_{pq}\in\bhpq$,
of~an \Mh function $f$ on~$\Bn$,
$$ \|f\|^2_\cici = \sum_{p,q} \frac{(p)_n(q)_n}{\Gamma(n)^2} \|f_{pq}\|^2_{\pBn}.  $$
\end{corollary}

\begin{definition} \label{PD}
We~call the space of \Mh functions $f=\sum_{p,q}f_{pq}$, $f_{pq}\in\bhpq$, on~$\Bn$ for which
$\|f\|^2_\cici 
 <+\infty$
the \emph{\Mh Dirichlet space}, denoted $\MM_\cici$.  \qed
\end{definition}

We~also denote by
$$ \MM_{\cici,0} := \{f\in\MM_\cici: f_{pq}=0 \text{ if } pq=0 \}  $$
the subspace of all functions in $\MM_\cici$ which have no \Ph component (i.e.~$Pf=0$).
The~quantity $\|f\|_\cici$ is a seminorm on $\MM_\cici$ and a norm on~$\MM_{\cici,0}$;
abusing the language, we~will often speak just of a ``norm''.
Obviously, $\MM_\cici$~contains all the spaces~$\bhpq$, $p,q\ge0$, and their span is dense in~it.

\begin{remark} The~space $\MM_{\cici,0}$ has reproducing kernel
$$ K^\cici(z,w) = \Gamma(n)^2 \sum_{p,q\ge1} \frac{K_{pq}(z,w)}{(p)_n(q)_n}.  $$
The~authors do not know if this sum can be evaluated explicitly.  \qed  \end{remark}

\section{Moebius invariance} \label{Sec4}
The~following facts likely are again quite standard, but we include the short proof
for completeness. Recall that for a function $f$ on~$\Bn$ and $0<r<1$, $f_r$~denotes
the function on $\pBn$ defined by $f_r(\zeta):=f(r\zeta)$.

\begin{lemma} \label{PE}
Let $f=\sum_{p,q}f_{pq}$, $f_{pq}\in\bhpq$, be an \Mh function on~$\Bn$.
Then the following three conditions are equivalent:
\begin{align}
& \sup_{0<r<1} \|f_r\|^2_\pBn <+\infty; \label{VH}  \\
& \text{a finite } \lim_{r\nearrow1} \|f_r\|^2_\pBn \text{ exists}; \label{VI} \\
& \sum_{p,q} \|f_{pq}\|^2_\pBn < +\infty .  \label{VJ}
\end{align}
The three quantities above are then equal, and the sum and the sequence
\begin{align}
& \sum_{p,q} f_{pq}  \label{VK}  \\
& \lim_{r\nearrow1} f_r  \label{VL}
\end{align}
then converge to the same function --- denoted~$f^*$ --- in~$L^2(\pBn,d\sigma)$;
$f^*$~can thus be interpreted as the ``boundary value'' of~$f$.
\end{lemma}

We~denote the quantity \eqref{VH}--\eqref{VJ} by $\|f\|^2_\Hardy$, and call the space of
all \Mh $f$ for which it is finite the \emph{\Mh Hardy space}~$\MM_\Hardy$. We~remark that
$\|f\|_\Hardy$ actually coincides with $\|f\|_{-1}$, and the reproducing kernel of $\MM_\Hardy$
was computed explicitly in~\cite{EY}.

Abusing the notation slightly, we~will sometimes write just~$f|_\pBn$, or~even~$f$, instead of~$f^*$,
and just $\|f\|_\pBn$ instead of $\|f^*\|_\pBn=\|f\|_\Hardy$.

\begin{proof} From \eqref{tVA} and~\eqref{tVG}, we~have
$$ \|f_r\|^2_\pBn = \sum_{pq} |S^{pq}(r)|^2 \|f_{pq}\|_\pBn^2 .  $$
Since $S^{pq}(r)$ are nondecreasing (strictly increasing for $pq>0$) functions of $r\in(0,1)$,
with $S^{pq}(1)=1$, it~follows that the limit \eqref{VI} coincides with the supremum~\eqref{VH};
and by the Lebesgue Monotone Convergence Theorem, they are both equal to~\eqref{VJ}.
This settles the first part.
For~the second, note that \eqref{VJ} means, by~\eqref{tVA}, that the partial sums of the series
\eqref{VK} form a Cauchy sequence; since $L^2$ is complete, they must have a limit~$f^*$,
and $\||f^*\|^2_\pBn=\sum_{p,q}\|f_{pq}\|^2_\pBn$. By~\eqref{tVA} and~\eqref{tVG} once again,
$$ \|f_r-f^*\|^2_\pBn = \sum_{p,q} (1-S^{pq}(r)^2) \|f_{pq}\|^2_\pBn .  $$
However, the right-hand side tends to zero again by the Lebesgue Monotone Convergence Theorem,
proving that $f_r\to f^*$.
\end{proof}

Recall that in the polar coordinates $z=r\zeta$ on~$\CC^n$ ($r>0$, $\zeta\in\pBn$),
the Euclidean Laplacian $\Delta$ is given~by
$$ \Delta = \frac{\partial^2}{\partial r^2} + \frac{2n-1}r \frac\partial{\partial r}
 + \frac1{r^2} \dsph,  $$
where $\dsph$ is the \emph{spherical Laplacian}, which involves only differentiations
with respect to the $\zeta$ variables. In~particular, the value of $\dsph\phi$ on
a sphere $|z|=$const. depends only on the values of the function $\phi$ on that sphere.
Another operator with this property is the \emph{complex normal derivative}
(or~\emph{Reeb vector field})
$$ \cR := \sum_{j=1}^n \Big(z_j\frac\partial{\partial z_j}-\oz_j\frac\partial{\partial\oz_j}\Big). $$
The operator $\dsph$ can be expressed explicitly~as
$$ \dsph = -\cR^2+ \sum_{j,k=1}^n (L_{jk}\overline L_{jk}+\overline L_{jk} L_{jk}) $$
where $L_{jk},\overline L_{jk}$ are the tangential vector fields
$$ L_{jk} := \oz_j\frac\partial{\partial z_k} - \oz_k\frac\partial{\partial z_j} , \qquad
 \overline L_{jk} := z_j\frac\partial{\partial\oz_k} - z_k\frac\partial{\partial\oz_j} .  $$
Both $\dsph$ and $\cR$ commute with the action of~$\Un$, i.e.~$\dsph(\phi\circ U)=(\dsph\phi)\circ U$
for any $U\in\Un$, and similarly for~$\cR$. (In~fact, the algebra of all $\Un$-invariant linear
differential operators on $\pBn$ is generated by $\dsph$ and $\cR$, but we will not need this fact.)
From the irreducibility of the multiplicity-free decomposition \eqref{tVB}, it~follows by abstract
theory that $\dsph$ and $\cR$ map each $\chpq$ (and~$\bhpq$) into itself and actually reduce on it to
a multiple of the identity. Evaluation on e.g.~the element $\zeta_1^p\overline\zeta_2^q\in\chpq$ (for~$n\ge2$)
shows that, explicitly,
\[ \begin{aligned}
\dsph|\chpq &= -(p+q)(p+q+2n-2) I|\chpq, \\
\cR|\chpq &= (p-q)I|\chpq
\end{aligned}  \label{tYA}   \]
(which prevail also for $n=1$ by the remarks on $\chpq$ when $n=1$ after~\eqref{tVG};
in~that case $\dsph=-\cR^2$).
Let~$\LL_m$, $1\le m\le2n(n-1)$, denote the collection of all the operators $L_{jk},\overline L_{jk}$,
$j,k=1,\dots,n$, $j\neq k$, in~some (fixed) order.

\begin{theorem} \label{PF}
If~$f$ is \Mh on $\Bn$, $n\ge2$, then $f\in\MM_\cici$ if and only if
\[ \sum_{j_1,j_2,\dots,j_n=1}^{2n(n-1)} \|\LL_{j_1}\LL_{j_2}\dots\LL_{j_n}(Qf)\|^2_\pBn <+\infty, \label{VM} \]
and the square root of the left-hand side is a seminorm equivalent to~$\|f\|_\cici$.
\end{theorem}

\begin{proof} Since $-\overline{L_{jk}}$ is the adjoint of $L_{jk}$ in $L^2(\pBn,d\sigma)$,
we~have for any $g\in C^2(\pBn)$
$$ \sum_{j=1}^{2n(n-1)} \|\LL_j g\|^2_\pBn = - \sum_{j,k=1}^n \spr{(L_{jk}\overline L_{jk}+\overline L_{jk} L_{jk})g,g}_\pBn
 = \spr{(-\dsph-\cR^2)g,g}_\pBn . $$
If~$g=\sum_{p,q}g_{pq}$, $g_{pq}\in\chpq$, is~the Peter-Weyl decomposition~\eqref{tVA} of~$g$,
we~thus have by~\eqref{tYA}
$$ \sum_{j=1}^{2n(n-1)} \|\LL_j g\|^2_\pBn = \sum_{p,q} [4pq+(2n-2)(p+q)] \|g_{pq}\|^2_\pBn.  $$
Iterating this formula, we~obtain
\[ \sum_{j_1,j_2,\dots,j_m=1}^{2n(n-1)} \|\LL_{j_1}\LL_{j_2}\dots\LL_{j_m}g\|^2_\pBn
 = \sum_{p,q} [4pq+(2n-2)(p+q)]^m \|g_{pq}\|^2_\pBn ,  \label{VN}  \]
for any $m=0,1,2,\dots$. Now for $p,q\ge1$ and $n\ge2$ obviously
\begin{align*}
[4pq+(2n-2)(p+q)]^n &\asymp [4pq+(2n-2)(p+q)+1]^n \\
&\asymp (p+1)^n(q+1)^n \\
&\asymp (p)_n (q)_n .
\end{align*}
Abusing notation by denoting by the same letter $Q$ also the orthogonal projection in
$L^2(\pBn,d\sigma)$ onto $\bigoplus_{p,q\ge1}\chpq$, we~thus have
$$ \sum_{j_1,j_2,\dots,j_n=1}^{2n(n-1)} \|\LL_{j_1}\LL_{j_2}\dots\LL_{j_n}(Qg)\|^2_\pBn
 \asymp \sum_{p,q=1}^\infty \frac{(p)_n(q)_n}{\Gamma(n)^2} \|g_{pq}\|^2_\pBn . $$
Applying this to $g=f_r$ and using the last lemma, the assertion follows.
\end{proof}

We~remark that, strictly speaking, in~\eqref{VM} we should write more correctly
\[ \sup_{0<r<1} \sum_{j_1,j_2,\dots,j_n=1}^{2n(n-1)} \|\LL_{j_1}\LL_{j_2}\dots\LL_{j_n}(Qf_r)\|^2_\pBn <+\infty, \]
as~$\LL_jf$ need not be harmonic in general when $f$~is; however, as $(\LL_jf)_r=\LL_j(f_r)$ (since $\LL_j$ are
tangential operators) and $(Qf)_r=Q(f_r)$ (from~\eqref{UY}), the~proof shows that the last supremum actually
coincides with $\|(-\dsph-\cR^2)^{n/2}f\|^2_\Hardy$, and the function $(-\dsph-\cR^2)^{n/2}f$
(defined in the sense of functional calculus for self-adjoint operators) is \Mh by~\eqref{tYA}.
We~thus take the liberty to use the shorthand indicated in the sentence before the proof of Lemma~\ref{PE}.

Yet another reformulation of the last theorem is as follows. Consider the weak-maximal operator $X$ acting
from $L^2(\pBn)$ into the Cartesian product of $2n^2(n-1)$ copies of $L^2(\pBn)$ by
$$ g\longmapsto \{\LL_{j_1}\dots\LL_{j_n}g\}_{j_1,\dots,j_n=1}^{2n(n-1)} ; $$
that~is, the~domain of $X$ consists of all $g\in L^2(\pBn)$ for which all the $\LL_{j_1}\dots\LL_{j_n}g$
exist in the sense of distributions and belong to $L^2(\pBn)$.
(In~other words, $X=Y^*$ where $Y$ is the restriction of the formal adjoint $X^\dagger$
of $X$ to $\bigoplus^{2n^2(n-1)} C^\infty(\pBn)$.)
Then~$f\in\MM$ belongs to $\MM_\cici$ if and only if $(Qf)^*\in\operatorname{Dom}(X)$, 
and $\|X(Qf)^*\|$ is a seminorm equivalent to~$\|f\|_\cici$.

\begin{corollary} \label{PG}
The space $\MM_\cici$ is Moebius invariant: $f\in\MM_\cici$ implies $f\circ\phi\in\MM_\cici$
for any $\phi\in\autBn$.
\end{corollary}

\begin{proof} Let $\phi\in\autBn$ and $f\in\MM_\cici$. We~need to show that the sum in \eqref{VM}
with $f\circ\phi$ in the place of $f$ is finite. Note that
$$ Q(f\circ\phi) = Q((Qf)\circ\phi)+Q((Pf)\circ\phi) = Q((Qf)\circ\phi) , $$
since composition with $\phi$ preserves holomorphy and, hence, plurisubharmonicity.
Note~further that, by~\eqref{VN}, for any $g\in L^2(\pBn,d\sigma)$,
\begin{multline*}
\sum_{j_1,j_2,\dots,j_m=1}^{2n(n-1)} \|\LL_{j_1}\LL_{j_2}\dots\LL_{j_m}Qg\|^2_\pBn
 = \sum_{p,q\ge1} [4pq+(2n-2)(p+q)]^m \|g_{pq}\|^2_\pBn   \\
\le \sum_{p,q\ge0} [4pq+(2n-2)(p+q)]^m \|g_{pq}\|^2_\pBn
 = \sum_{j_1,j_2,\dots,j_m=1}^{2n(n-1)} \|\LL_{j_1}\LL_{j_2}\dots\LL_{j_m}g\|^2_\pBn;
\end{multline*}
hence it is enough to show that, in~fact, even the sum in \eqref{VM} with $(Qf)\circ\phi$
in the place of $Qf$ is finite.

Observe that the tangential vector-fields $\LL_m$, $m=1,\dots,2n(n-1)$, span (very redundantly)
the~entire complex tangent space to~$\pbn$. Thus for any differentiable function $g$ on~$\pBn$,
$\sum_m\|\LL_m g\|^2_\pBn\asymp\|\nct g\|^2$, the~norm-square of the restriction $\nct g$ of
the tensor $\nabla g$ to the complex tangent space of~$\pBn$ in the sense of complex geometry.
Now~for any vector field $X$ on~$\pBn$, one~has $X(g\circ\phi)=d\phi(X)g$.
Since $\phi$ maps the sphere $\pBn$ onto itself, the~derived map $d\phi$
maps the real tangent space of $\pBn$ into itself. As~$\phi$ is holomorphic,
$d\phi$~is complex linear, hence $d\phi$ maps also the complex tangent space
(consisting of all vectors $Y$ such that both $Y$ and $iY$ are real tangent)
of~$\pBn$ into itself. Finally, $d\phi|\pBn$ is a smoothly varying map on the
compact manifold~$\pBn$ (hence, in~particular, so~is its Jacobian).
Consequently,
$$ \|\nct(g\circ\phi)\|^2_\pBn = \|d\phi(\nct) g\|^2_\pBn \le C_\phi \|\nct g\|^2_\pBn  $$
with some finite $C$ (independent of~$g$). Iterating this argument, it~transpires that
$$ \|\nct^n(g\circ\phi)\|^2_\pBn \le C_\phi^n \|\nct g\|^2_\pBn.  $$
Passing from $\nct^n$ back to the~$\LL_m$, the~last inequality reads
$$ \sum_{j_1,j_2,\dots,j_n=1}^{2n(n-1)} \|\LL_{j_1}\LL_{j_2}\dots\LL_{j_n}(g\circ\phi)\|^2_\pBn
 \le C_\phi^n \sum_{j_1,j_2,\dots,j_n=1}^{2n(n-1)} \|\LL_{j_1}\LL_{j_2}\dots\LL_{j_n}g\|^2_\pBn.  $$
Taking $g=Qf$, the proof is complete.
\end{proof}

\begin{remark} The~authors suspect that, for any fixed $s>-n-1$,
\[ C_{pq}(s) \asymp [(p+1)(q+1)]^{-s-1}  \label{VO}  \]
uniformly for all $p,q\in\NN$. If~this is rue, then the proof of Theorem~\ref{PF} shows that
the condition \eqref{VM} is further equivalent~to
$$ \sum_{j_1,j_2,\dots,j_{n+k+1}=1}^{2n(n-1)} \|\LL_{j_1}\LL_{j_2}\dots\LL_{j_{n+k+1}}(Qf)\|^2_k <+\infty $$
for some (equivalently, any) nonnegative integer~$k$. See~the proof of Theorem~\ref{PH} below for the details.
The~authors showed in Theorem~11 in~\cite{EY} that \eqref{VO} holds when $p,q$ tend to infinity with
he ratio $p/q$ fixed, but were unable to get a uniform estimate.  \qed  \end{remark}

\begin{remark} Another consequence of \eqref{VO} would be an extension of the definition of
$\MM_s$ from the original range $s>-1$, and our ``analytic continuation'' to $s>-n-1$,
to~all real~$s$. Namely, denoting
$$ \MM_{\#s} :=\{f=\sum_{p,q}f_{pq}, \; f_{pq}\in\bhpq: \; \sum_{p,q} \frac{\|f_{pq}\|^2_\pBn}{(p+1)^{s+1}(q+1)^{s+1}} =: \|f\|^2_{\#s}<+\infty\}  $$
we~would then have
$$ \MM_{\#s} = \MM_s \quad\text{for $s>-n-1$, with equivalent norms,}  $$
by~\eqref{VO}, and
$$ Q\MM_{\#s} = Q\MM_\cici \quad\text{for $s=-n-1$,}  $$
with equivalent norms on~$M_{\cici,0}$, by~\eqref{VG}. Since evidently $\MM_{\#s}\subset\MM_{\#s'}$ continuously for $s<s'$,
it~would also follow that
$$ \MM_{\cici,0} \subset \MM_s \qquad \forall s>-n-1,  $$
which inclusion the current authors are unable to verify.

(We~believe the lower bound in \eqref{VO} can be obtained from the inequality
$$ \frac{\Gamma(n+p)\Gamma(n+q)}{\Gamma(n)\Gamma(n+p+q)} \FF21{p,q}{p+q+n}t \ge \frac{t^{pq}}{3^{1/4}}, $$
which seems to be true but we have not been able to prove~it. We~also have no clue how to get the upper bound in~\eqref{VO}.)  \qed  \end{remark}

A~priori, it~is not evident that
$$ \|f\|^2_\cici = \sum_{P,q} \|f_{pq}\|^2_\pBn \lim_{s\searrow-n-1} (n+s+1)^2 C_{pq}(s)  $$
coincides with
$$  \lim_{s\searrow-n-1} (n+s+1)^2 \sum_{P,q} C_{pq}(s) \|f_{pq}\|^2_\pBn =  \lim_{s\searrow-n-1} (n+s+1)^2 \|f\|^2_s  $$
--- none of the standard conditions for interchanging the limit and the summation seems to apply.
However, by~Fatou's lemma, we~at least always have
\[ \|f\|^2_\cici \le \liminf_{s\searrow-n-1} (n+s+1)^2 \|f\|^2_s , \label{VX}  \]
with equality, of~course, when the sums above are finite (i.e.~for $f$ with only finitely
many nonzero $\bhpq$-components). In~other words, if~we introduce the space
\[ \MM' := \{ f\in\MM: f\in\MM_s\;\forall s>-n-1\text{ and a finite} \lim_{s\searrow-n-1} (n+s+1)^2 \|f\|^2_s \text{ exists}\},  \label{XA} \]
then plainly the last limit is a (semi-)norm on~$\MM'$, the~limit
$$ \spr{f,g}' := \lim_{s\searrow-n-1} (n+s+1)^2 \spr{f,g}_s  $$
exists for any $f,g\in\MM'$ and makes $\MM'$ into a (semi-)inner product space, and by the remarks above,
\[ \MM'\subset\MM_\cici \text{ continuously}  \label{XB}  \]
while the algebraic span of~$\bhpq$, $p,q\ge0$, is~contained in $\MM'$ and the norms
$\|\cdot\|'$ and $\|\cdot\|_\cici$ coincide on~it. It~follows that $\MM_\cici$ is just
the completion of $\MM'$ with respect to the above norm.

\begin{remark} It~follows from \eqref{XC} in Remark~\ref{XXCC} above that $\|f\|^2_s$
is actually a non-increasing function of $s>-n-1$, for any~$f\in\MM$. In~particular, one~has
$$ \MM_s \subset \MM_{s'} \quad\text{if } s'>s>-n-1.   \qed  $$
\end{remark}

\begin{remark} Up~to the authors' knowledge, it seems to be an interesting open problem
whether $\lim_{s\searrow-n-1} (n+s+1)^2 \|f\|^2_s$ actually always exists for $f\in\MM_\cici$
and coincides with~$\|f\|^2_\cici$ (that~is, whether $\MM'=\MM_\cici$).
The~analogous assertion in the context of the ordinary holomorphic Dirichlet space holds:
namely, by~\eqref{TD},
$$ (n+s+1)\|f\|^2_s = (n+s+1) \|f_0\|^2_\pBn + \sum_{j\ge1} \frac{(n)_j}{(n+s+2)_{j-1}} \|f_j\|^2_\pBn ,  $$
and as $s\searrow-n-1$, the last sum tends to $\sum_{j\ge1} \frac{(n)_j}{(j-1)!} \|f_j\|^2_\pBn = \|f\|^2_\circ$
by the Lebesgue Monotone Convergence Theorem.   \qed  \end{remark}

\begin{theorem} \label{PK}
The inner product in $\MM_\cici$ is Moebius invariant:
\[ \spr{f,g}_\cici = \spr{f\circ\phi,g\circ\phi}_\cici   \label{XD}  \]
for any $f,g\in\MM_\cici$ and $\phi\in\autBn$.
\end{theorem}

\begin{proof} Any $\phi\in\autBn$ can be written in the form $\phi=U\circ\phi_a\circ V$,
where $U,V\in\Un$ while $\phi_a$ denotes the geodesic symmetry (i.e.~$\phi_a\circ\phi_a=\text{id}$
and $\phi_a$ has only an isolated fixed-point) interchanging the origin $0\in\Bn$ with the point
$(a,0,0,\dots,0)\in\Bn$, $0\le a<1$.
Since both $\MM_\cici$ and its inner product $\spr{\cdot,\cdot}_\cici$ are $\Un$-invariant
(by~their very construction), it~is enough to prove the assertion for $\phi=\phi_a$.
Furthermore, since we know from Corollary~\ref{PG} (or,~rather, from its~proof)
that the composition operator $f\mapsto f\circ\phi$ is continuous on~$\MM_\cici$,
it~is further enough to prove the assertion for $f,g$ in a dense subset of~$\MM_\cici$.
In~particular, by~linearity, we~may assume that $f\in\bhpq$ and $g\in\mathbf H^{p'q'}$
for some $p,q,p',q'\in\NN$.

We~will show that under all these hypotheses, $\spr{f\circ\phi_a,g\circ\phi_a}'$ exists
for all $0\le a<1$ and does not depend on~$a$. By~the observations in the paragraph
before the theorem, this will complete the proof.

Fix $0<\rho<1$. Recall that the measure
$$ d\tau(z) := \frac{dz}{(1-|z|^2)^{n+1}}  $$
on~$\Bn$ is invariant under~$\phi$, and also
$$ 1-|\phi_a(z)|^2 = \frac{(1-|a|^2)(1-|z|^2)}{|1-az_1|^2}.  $$
By~the change of variable $z\mapsto\phi_a(z)$, we~thus have, for any $s>-1$,
\begin{align*}
\spr{f\circ\phi_a,g\circ\phi_a}_s
&= \frac{(s+1)_n}{\pi^n} \int_\Bn (f\overline g)(\phi_a(z)) (1-|z|^2)^{s+n+1}\,d\tau(z) \\
&= \frac{(s+1)_n}{\pi^n} \int_\Bn (f\overline g)(z) (1-|\phi_a(z)|^2)^{s+n+1}\,d\tau(z) \\
&= \int_\Bn (f\overline g)(z) \Big(\frac{1-a^2}{|1-az_1|^2}\Big)^{n+s+1} \,d\mu_s(z).
\end{align*}
Passing to the polar coordinate $z=r\zeta$, with $0\le r<1$ and $\zeta\in\pBn$, we~can continue with
\[ = \frac{(s+1)_n}{\pi^n} \int_0^1 \frac{2\pi^n}{\Gamma(n)} \int_\pBn (f\overline g)(r\zeta)
 \Big(\frac{1-a^2}{|1-ar\zeta_1|^2}\Big)^{n+s+1} (1-r^2)^s r^{2n-1} \,d\zeta \,dr, \label{XE} \]
that~is, using~\eqref{tVG}.
\begin{align*}
&= \frac{(s+1)_n}{\pi^n} \int_0^1 G(a,r) (1-r^2)^s r^{2n-1} \,dr , \\
& \hskip4em \text{where } G(a,r) := \frac{2\pi^n}{\Gamma(n)} S^{pq}(r) S^{p'q'}(r)
 \int_\pBn (f\overline g)(\zeta) \Big(\frac{1-a^2}{|1-ar\zeta_1|^2}\Big)^{n+s+1} \,d\zeta .
\end{align*}
Carrying out the $\zeta$ integration shows that $G(a,r)$ is a holomorphic function of $|a|<\rho$
and $|r|<1/\rho$. Invoking Lemma~\ref{PB}, it~thus follows in the same way as in the proof of
Proposition~\ref{PC} that $\spr{f\circ\phi_a,g\circ\phi_a}_s$
extends to a holomorphic function of $|a|<\rho$ and $s\in\CC$, except for at most double poles at
$s=-n-1,-n-2,\dots,-2n$ and at most triple poles at $s=-2n-j-1$, $j\in\NN$.
Consequently, the~function $(s+n+1)^2\spr{f\circ\phi_a,g\circ\phi_a}_s$ extends to a holomorphic
function of $|a|<\rho$ and $s\in\CC$ except for poles as above, excluding $s=-n-1$ where it
assumes a finite value. In~particular (taking $f=g$), this means that $f\circ\phi_a,g\circ\phi_a\in\MM'$
for all $0\le a<\rho$, and the inner product $\spr{f\circ\phi_a,g\circ\phi_a}'$ is a smooth
function of these~$a$.

Finally, it~is legitimate to differentiate under the integral sign in~\eqref{XE}, yielding,
for $s>-1$,
\[ \begin{aligned}
& \frac\partial{\partial a} \spr{f\circ\phi_a,g\circ\phi_a}_s = \frac{(s+1)_n}{\pi^n}
 \int_0^1 \frac{2\pi^n}{\Gamma(n)} \int_\pBn (f\overline g)(\zeta) S^{pq}(r) S^{p'q'}(r) \;\times \\
& \hskip4em (n+s+1) \Big(\frac{1-a^2}{|1-ar\zeta_1|^2}\Big)^{n+s} \Big[\frac\partial{\partial a}
 \frac{1-a^2}{|1-ar\zeta_1|^2} \Big] (1-r^2)^s r^{2n-1} \,d\zeta \,dr . \end{aligned}  \label{XF}  \]
Repeating the argument above, it~transpires that for all $0\le a<\rho$,
\[ \frac\partial{\partial a} (n+s+1)^2 \spr{f\circ\phi_a,g\circ\phi_a}_s = (n+s+1) F_a(s),  \label{XG} \]
where $F_a(s)$ is a holomorphic function of $s$ except for at most triple poles at $s=-2n-1-j$, $j\in\NN$,
and at most double poles at $s=-n-2,\dots,-2n$; in~particular, $F_a(s)$~is holomorphic near $s=-n-1$
and assumes a finite value there. Hence, thanks to the factor $n+s+1$ in~\eqref{XG},
$$ \frac\partial{\partial a} \spr{f\circ\phi_a,g\circ\phi_a}' =0 \quad\text{for }0\le a<\rho.  $$
Since $\rho$ was arbitrary, it~follows that $\spr{f\circ\phi_a,g\circ\phi_a}'=
\spr{f\circ\phi_0,g\circ\phi_0}'=\spr{f,g}'$ for all $0\le a<1$, completing the proof.
\end{proof}

\section{The \Ph Dirichlet space}  \label{Sec5}
For $pq=0$, the coefficients $C_{p0}(s)=C_{0p}(s)$ have only a single pole at $s=-n-1$
(cf.~\eqref{VP}) for $p\neq0$, with residue
$$ C^\circ_{p0} := \lim_{s\searrow-n-1}(n+s+1)C_{p0}(s) = \frac{(n)_p}{\Gamma(p)},  $$
while $C^\circ_{00} := \lim_{s\searrow-n-1}(n+s+1)C_{00}(s) =0$.
Accordingly, $\MM_\circ$~consists only of \Ph functions, with (semi-)norm
$$ \|f\|^2_\circ := \sum_{p=1}^\infty p\frac{(n)_p}{p!} (\|f_{p0}\|^2_\pBn+\|f_{0p}\|^2_\pBn). $$
In~other words,
$$ \MM_\circ = \AA_\circ \oplus \overline{\AA_\circ}  $$
is just the orthogonal sum of the usual holomorphic Dirichlet space $\AA_\circ$ and its complex conjugate.

The~result below, parallel to Theorem~\ref{PF} for the \Mh case, is~likely folk lore,
but~the authors are unaware of a specific reference.

\begin{theorem} \label{PH}
If $f$ is \Ph on $\Bn$, $n\ge2$, then $f\in\MM_\circ$ if and only if
\[ \sum_{j_1,j_2,\dots,j_n=1}^{2n(n-1)} \|\LL_{j_1}\LL_{j_2}\dots\LL_{j_n}f\|^2_\Hardy <+\infty \label{VQ} \]
if and only if
\[ \sum_{j_1,j_2,\dots,j_{n+k+1}=1}^{2n(n-1)} \|\LL_{j_1}\LL_{j_2}\dots\LL_{j_{n+k+1}}f\|^2_k <+\infty, \label{VR}  \]
for some (equivalently, any) nonnegative integer~$k$.
\end{theorem}

\begin{proof} As~we have seen in \eqref{VN} in the proof of Theorem~\ref{PF}, \eqref{VQ}~equals,
for $f=\sum_{p,q}f_{pq}$ with $f_{pq}\in\bhpq$,
$$ \sum_{j_1,j_2,\dots,j_n=1}^{2n(n-1)} \|\LL_{j_1}\LL_{j_2}\dots\LL_{j_n}f\|^2_\pBn
 = \sum_{p,q} [4pq+(2n-2)(p+q)]^n \|f_{pq}\|^2_\pBn . $$
Since $f$ is now pluriharmonic, the right-hand side reduces just~to
$$ \sum_{p=1}^\infty [(2n-2)p]^n (\|f_{p0}\|^2_\pBn+\|f_{0p}\|^2_\pBn) .  $$
As $n\ge2$ by hypothesis, we~have
$$ [(2n-2)p]^n \asymp p^n \asymp \frac{\Gamma(n+p)}{\Gamma(n)\Gamma(p)} = C^\circ_{p0} $$
for all~$p$, and the first claim follows.

For~the second claim, denote again, for any function $f$ on~$\Bn$, $f_r(\zeta):=f(r\zeta)$
for $0<r<1$ and $\zeta\in\pBn$. Using \eqref{VN} for $g=f_r$ yields
$$ \sum_{j_1,j_2,\dots,j_m=1}^{2n(n-1)} \|\LL_{j_1}\LL_{j_2}\dots\LL_{j_m}f_r\|^2_\pBn
 = \sum_{p,q} [4pq+(2n-2)(p+q)]^m S^{pq}(r)^2 \|f_{pq}\|^2_\pBn , $$
since $(f_r)_{pq}(\zeta)=S^{pq}(r)f_{pq}(\zeta)$ by~\eqref{tVG}.
As~$\LL_j$, being tangential, do~not act on the $r$ variable, we~also have
$$ \LL_{j_1}\LL_{j_2}\dots\LL_{j_m}f_r = (\LL_{j_1}\LL_{j_2}\dots\LL_{j_m}f)_r . $$
Hence for any $s>-1$,
\begin{align*}
\|\LL_{j_1}\LL_{j_2}\dots\LL_{j_m}f\|^2_s
&= \frac{(s+1)_n}{\pi^n} \int_0^1 \frac{2\pi^n}{\Gamma(n)} \|(\LL_{j_1}\LL_{j_2}\dots\LL_{j_m}f)_r\|^2_\pBn (1-r^2)^s r^{2n-1}\,dr \\
&= \frac{(s+1)_n}{\Gamma(n)} \sum_{p,q} [4pq+(2n-2)(p+q)]^m \|f_{pq}\|^2_\pBn \int_0^1 S^{pq}(\sqrt t)^2 t^{n-1} (1-t)^s \,dt \\
&= \sum_{p,q} [4pq+(2n-2)(p+q)]^m \|f_{pq}\|^2_\pBn C_{pq}(s).
\end{align*}
Specializing now to the current \Ph case $pq=0$, we~again have for all $p\ge1$
$$ [(2n-2)p]^m C_{p0}(s) \asymp p^m C_{p0}(s) \asymp p^{m-s-1} \asymp p^{m-n-s-1}C^\circ_{p0} . $$
Hence for $s=k$ and $m=n+k+1$, with any $k=0,1,2,\dots$,
$$ [(2n-2)p]^{n+k+1} C_{p0}(k) \asymp C^\circ_{p0}  \qquad\forall p\in\NN  $$
(for $p=0$, both sides vanish), and the second claim follows.
\end{proof}

The~following simple result seems to have no counterpart in the \Mh case.

\begin{theorem} \label{PI}
If $f$ is \Ph on $\Bn$, $n\ge1$, then $f\in\MM_\circ$ if and only if
$$ \|\cN^m f\|^2_{2m-n-1} <+\infty  $$
for some (equivalently, any) integer~$m>\frac n2$.
Furthermore, the square root of the left-hand side is a seminorm equivalent to~$\|f\|_\circ$.
\end{theorem}

\begin{proof} By~straightforward inspection,
$$ \cN f_{p0} = pf_{p0}, \quad \cN f_{0p}=pf_{0p} \qquad\forall p\ge0.  $$
Consequently,
$$ \|\cN^m f_{pq}\|^2_\pBn = (p+q)^{2m} \|f_{pq}\|^2_\pBn \qquad\text{for }pq=0,  $$
and, as~in the preceding proof, for any \Ph $f$,
\begin{align*}
\|\cN^m f\|^2_s &= \sum_p p^{2m} C_{p0}(s) (\|f_{p0}\|^2_\pBn+\|f_{0p}\|^2_\pBn) \\
&\asymp p^{2m-s-n-1} C^\circ_{p0} (\|f_{p0}\|^2_\pBn+\|f_{0p}\|^2_\pBn) \\
&\asymp \|f\|_\circ^2 \quad\text{if }2m=n+s+1,
\end{align*}
completing the proof.
\end{proof}

\section{The harmonic Dirichlet space} \label{Sec6}
The~situation in the harmonic case is pretty similar as for the \Ph functions
in the preceding section, so~we will be brief.
For~all $p\ge0$, let $\bhp$ be the space of harmonic polynomials on~$\RR^n$, $n\ge2$,
homogeneous of degree~$p$, and let $\chp$ be the space of restrictions of elements
of $\bhp$ to the unit sphere~$\pbn$. We~refer to \cite{ABR}, especially Chapter~5,
for the Peter-Weyl decomposition
\[ L^2(\pbn,d\sigma) = \bigoplus_{p=0}^\infty \chp   \label{WA}  \]
under the action of the orthogonal group $O(n)$ of rotations of~$\RR^n$,
and the associated decomposition
\[ \HH = \bigoplus_p \bhp   \label{WB}  \]
of the space of all harmonic functions on the unit ball $\bn$ of $\RR^n$ into the
direct sum of the~$\bhp$: namely, any harmonic function $f$ on $\bn$ can be
uniquely written~as
\[ f = \sum_{p=0}^\infty f_p, \qquad f_p\in\bhp,  \label{WC} \]
with the sum converging uniformly on compact subsets. Here $d\sigma$ now stands for
the normalized surface measure on~$\pbn$. The~\emph{weighted harmonic Bergman space}
$$ \HH_s(\bn) := \{ f\in L^2(\bn,d\rho_s): \;f\text{ is harmonic on }\bn\}  $$
consists of all harmonic functions on $\bn$ square-integrable with respect to the measure
\[ d\rho_s(x) := \frac{\Gamma(\frac n2+s+1)}{\pi^{n/2}\Gamma(s+1)} (1-|x|^2)^s\,dx,  \qquad s>-1, \label{WD} \]
where $dx$ denotes the Lebesgue volume on~$\RR^n$. The~restriction on $s$ ensures that these spaces are nontrivial,
and the factor $\frac{\Gamma(\frac n2+s+1)}{\pi^{n/2}\Gamma(s+1)}$ makes $d\rho_s$ a probability measure,
so~that $\|\jedna\|=1$. For~$f$ as in~\eqref{WC}, we~have by the orthogonality in~\eqref{WA}
\begin{align}
\|f\|^2_s &= \frac{\Gamma(\frac n2+s+1)}{\pi^{n/2}\Gamma(s+1)} \int_0^1 \frac{2\pi^{n/2}}{\Gamma(\frac n2)} \int_\pbn
 |f(r\zeta)|^2 \,d\sigma(\zeta) \,(1-r^2)^s r^{n-1}\,dr \nonumber \\
&= \frac{\Gamma(\frac n2+s+1)}{\Gamma(\frac n2)\Gamma(s+1)} \int_0^1 \sum_{p=0}^\infty r^{2p} \|f_p\|^2_\pbn (1-r^2)^s \,2r^{n-1}\,dr \nonumber \\
&= \frac{\Gamma(\frac n2+s+1)}{\Gamma(\frac n2)\Gamma(s+1)} \sum_p \|f_p\|^2_\pbn \int_0^1 t^{p+\frac n2-1} (1-t)^s \,dt \nonumber \\
&= \frac{\Gamma(\frac n2+s+1)}{\Gamma(\frac n2)\Gamma(s+1)} \sum_p \frac{\Gamma(p+\frac n2)\Gamma(s+1)}{\Gamma(p+s+\frac n2+1)} \|f_p\|^2_\pbn \nonumber \\
&= \sum_p \frac{(\frac n2)_p}{(\frac n2+s+1)_p} \|f_p\|^2_\pbn ,  \label{WE}
\end{align}
and, accordingly, the reproducing kernel of $\HH_s$ is given~by
\[ K^\harm_s(x,y) = \sum_p \frac{(\frac n2+s+1)_p}{(\frac n2)_p} Z_p(x,y),  \label{WF} \]
where $Z_p(x,y)$, the reproducing kernel of~$\bhp$, is given by so-called zonal harmonics
(expressible explicitly in terms of Gegenbauer polynomials); see Chapter~8 in~\cite{ABR}
for the unweighted case, the~weighted case being completely parallel.

The~coefficients
\[ C^\harm_p(s) := \frac{(\frac n2)_p}{(\frac n2+s+1)_p}, \qquad p\in\NN,  \label{WG}  \]
extend to nonvanishing holomorphic functions of $s$ on the entire~$\CC$, except for simple poles
at $s=-\frac n2-1,\dots,-\frac n2-p$; accordingly, $K^\harm_s(x,y)$ extends to a holomorphic
function of~$s\in\CC$. Due~to the orthogonality of the spaces~$\bhp$, the~extended kernel ---
still denoted~$K^\harm_s$ --- will remain positive definite as long as $1/C^\harm_p(s)\ge0$
$\forall p\in\NN$, hence, precisely for $s\in\itv{-\frac n2-1,+\infty}$. The~last interval
is thus the ``harmonic Wallach set'' of~$\bn$.
The~norm in the corresponding reproducing kernel Hilbert spaces --- still denoted by~$\HH_s$ ---
is~still given by \eqref{WE} for $s>-\frac n2-1$. For $s=-\frac n2-1$, \eqref{WF}~reduces to
constant~one, and the associated space thus consists only of the constants, with $\|\jedna\|=1$.
As~the ``residue'' at $s=-\frac n2-1$, we~get the reproducing kernel
\[ K^\harm_\sq (x,y) := \lim_{s\searrow-\frac n2-1} \frac{K^\harm_s(x,y)-1}{s+\frac n2+1} . \label{WH}  \]
The~corresponding reproducing kernel Hilbert space $\HH_\sq$ consists of all harmonic functions
on $\bn$ for which
\[ \|f\|^2_\sq := \sum_{p=0}^\infty C^\sq_p \|f_p\|^2_\pbn < +\infty,  \label{WI} \]
where
\[ C^\sq_p := \lim_{s\searrow-\frac n2-1} (s+\tfrac n2+1) C^\harm_p(s) = p\frac{(\frac n2)_p}{p!}. \label{WJ}  \]
This can be viewed as the \emph{harmonic Dirichlet space}. It~is easily seen to coincide with
the eponymous space studied by other authors, see~e.g.~\cite{GKU} and the numerous references therein.
The~characterization of $\HH_\sq$ given in Theorem~\ref{PJ} below, however, seems not to appear
in the literature (up~to the authors' knowledge).

\begin{remark} As~in the holomorphic case, the~limit $\lim_{s\searrow-\frac n2-1}\|f\|^2_s$ always
exists for any $f\in\HH_\sq$ and coincides with~$\|f\|^2_\sq$. The~proof is the same as for the
holomorphic case.  \qed  \end{remark}

Finally, the~following characterization of $\HH_\sq$ can be given along the same lines as in
the preceding sections. Recall our notation
$$ X_{jk} = x_j\partial_k - x_k\partial_j , \qquad j,k=1,\dots,n, \quad j\neq k,  $$
for the tangential vector fields on~$\RR^n$, and $\XX_m$, $m=1,\dots,n(n-1)$, for~the collection
of all the~$X_{jk}$ (in~some fixed order). By~a~routine computation, one~checks that
$$ \sum_{j,k=1}^n X_{jk}^2 = 2\dsph,  $$
where $\dsph$ is the spherical Laplacian on~$\RR^n$: for $x=r\zeta$ with $r>0$ and $\zeta\in\pbn$,
$$ \Delta = \frac{\partial^2}{\partial r^2} + \frac{n-1}r \frac\partial{\partial r} + \frac1{r^2}\dsph.  $$
The~operator $\dsph$ commutes with the action of the orthogonal group $O(n)$ of~$\RR^n$,
hence it is automatically diagonalized by the Peter-Weyl decomposition~\eqref{WA}:
a~simple  computation reveals that
\[ \dsph|\bhp = -p(p+n-2)I|\bhp  \label{WK}  \]
where $I$ stands for the identity operator.

\begin{theorem} \label{PJ}
If $f$ is harmonic on~$\bn$, $n\ge2$, then $f\in\HH_\sq$ if and only~if
$$ \sum_{j_1,\dots,j_m=1}^{n(n-1)} \|\XX_{j_1}\dots\XX_{j_m}f\|^2_{2m-\frac n2-1} <+\infty . $$
for some (equivalently, any) integer $m>\frac n4$.
Furthermore, the square roots of the left-hand sides are seminorms equivalent to~$\|f\|_\sq$.
\end{theorem}

\begin{proof} Since the adjoint of $X_{jk}$ in $L^2(\pbn,d\sigma)$ is just~$-X_{jk}$,
we~have for any $g\in L^2(\pbn,d\sigma)$
$$ \sum_{j=1}^{n(n-1)} \|\XX_j g\|^2_\pbn = -\sum_{j,k=1}^n \spr{X_{jk}^2g,g}_\pbn = -2\spr{\dsph g,g}_\pbn, $$
so~for $g=\sum_p g_p$, $g_p\in\bhp$, as~in~\eqref{WA},
\[ \sum_{j=1}^{n(n-1)} \|\XX_j g\|^2_\pbn = \sum_p 2p(p+n-2) \|g_p\|^2_\pbn ,  \label{WL} \]
by~\eqref{WK}. Iterating this procedure, we~get
$$ \sum_{j_1,\dots,j_m=1}^{n(n-1)} \|\XX_{j_1}\dots\XX_{j_m}g\|^2_\pbn = \sum_p [2p(p+n-2)]^m \|g_p\|^2_\pbn . $$
Applying this now to $g(\zeta)=f(r\zeta)$ where $f$ is harmonic on~$\bn$, we~obtain
$$ \sum_{j_1,\dots,j_m=1}^{n(n-1)} \|\XX_{j_1}\dots\XX_{j_m}f(r\cdot)\|^2_\pbn = \sum_p [2p(p+n-2)]^m r^{2p} \|f_p\|^2_\pbn , $$
and, as~in~\eqref{WE}, for any $s>-1$,
\begin{align*}
& \sum_{j_1,\dots,j_m=1}^{n(n-1)} \|\XX_{j_1}\dots\XX_{j_m}f\|^2_s \\
& \hskip4em = \frac{\Gamma(\frac n2+s+1)}{\pi^{n/2}\Gamma(s+1)} \int_0^1 \frac{2\pi^{n/2}}{\Gamma(\frac n2)}
 \sum_{j_1,\dots,j_m=1}^{n(n-1)} \|\XX_{j_1}\dots\XX_{j_m}f(r\cdot)\|^2_\pbn \,(1-r^2)^s r^{n-1}\,dr \\
& \hskip4em = \frac{\Gamma(\frac n2+s+1)}{\Gamma(\frac n2)\Gamma(s+1)} \int_0^1 \sum_p [2p(p+n-2)]^m \|f_p\|^2_\pbn t^{p+\frac n2-1}(1-t)^s \,dt \\
& \hskip4em = \sum_p [2p(p+n-2)]^m \frac{(\frac n2)_p}{(\frac n2+s+1)_p} \|f_p\|^2_\pbn .
\end{align*}
Now for all $p\ge1$
$$ [2p(p+n-2)]^m \frac{(\frac n2)_p}{(\frac n2+s+1)_p} \asymp p^{2m-s-1} \asymp p^{2m-s-1-\frac n2} C^\sq_p , $$
whence
$$ [2p(p+n-2)]^m \frac{(\frac n2)_p}{(\frac n2+s+1)_p} \asymp C^\sq_p $$
if $2m=s+1+\frac n2$ (for $p=0$, both sides vanish). This completes the proof. 
\end{proof}

\end{document}